\documentclass[12pt,a4paper]{article}
\usepackage{amsthm}
\usepackage{fleqn}
\usepackage{array}
\usepackage[latin1]{inputenc}
\usepackage{amssymb}
\usepackage{amsfonts}
\usepackage{amsmath}

\begin{document}

\renewcommand{\thefootnote}{\fnsymbol{footnote}}
\newtheorem{thm}{Theorem}[section]
\newtheorem{prop}{Proposition}[section]
\newtheorem{lem}{Lemma}[section]
\newtheorem{cor}{Corollary}[section]

\begin{center}
{\large \bf {Curvature types of planar curves for gauges}}

\vspace{5mm}

Vitor Balestro, Horst Martini and Makoto Sakaki\footnote{Corresponding author}
\end{center}

\vspace{3mm}

{\small
Vitor Balestro, Instituto de Matematica e Estatistica, Universidade Federal Fluminense, 24210201 Niteroi, Brazil

e-mail: vitorbalestro@id.uff.br

\vspace{3mm}

Horst Martini, Fakult$\ddot{\mbox{a}}$t f$\ddot{\mbox{u}}$r Mathematik, Technische Universit$\ddot{\mbox{a}}$t Chemnitz, 09107 Chemnitz, Germany

e-mail: martini@mathematik.tu-chemnitz.de

\vspace{3mm}

Makoto Sakaki, Graduate School of Science and Technology, Hirosaki University, Hirosaki 036-8561, Japan

e-mail: sakaki@hirosaki-u.ac.jp
}

\vspace{5mm}

\textbf{Abstract.} In this paper results from the differential geometry of curves are 
extended from normed planes to gauge planes which are obtained by 
neglecting the symmetry axiom. Based on the gauge analogue of the 
notion of Birkhoff orthogonality from Banach space theory, we study 
all curvature types of curves in gauge planes, thus
generalizing their complete classification for normed planes. We 
show that (as in the subcase of normed planes) there are four such types, 
and we call them analogously Minkowski, normal, circular, and 
arc-length curvature. We study relations between them and extend, based on this, 
also the notions of evolutes and involutes to gauge planes.

\vspace{3mm}

\textbf{Mathematics Subject Classification (2010).}  46B20, 52A10,  
52A21, 53A04, 53A35

\vspace{3mm}

\textbf{Keywords.} Birkhoff orthogonality, curvature, 
evolute and involute, gauge, generalized Minkowski plane, planar curve

\section{Introduction}

Although the concept of normed (or real Banach) planes describes the  
setting of Finsler geometry locally, the knowledge on existing  
types of curvatures of curves in such planes is widespread and not  
systematized in the literature; only very recently this gap was filled  
by the paper \cite{BMS}. Using the notion of Birkhoff orthogonality from  
Banach space theory in a meaningful way, the authors of \cite{BMS} gave the  
first complete representation of all four curvature types making sense for 
curves in such planes. They also studied geometric properties of these 
types and their
relations to each other. In the present paper, we extend  
this framework from normed to gauge planes, i.e., to planes satisfying all the  
axioms of two-dimensional real Banach spaces except for the symmetry  
axiom. For this purpose, we introduce an analogously generalized  
orthogonality type (as for norms, we stay with the name "Birkhoff  
orthogonality"), and the same is done with other needed notions. It  
turns out that we obtain precisely four natural analogues of the curvature types  
for norms, and correspondingly we call them Minkowski, normal,  
circular, and arc-length curvature. We present basic geometric relations  
between these curvature types, and as an application the concepts of  
evolutes and involutes as well as reverse evolutes and reverse  
involutes are introduced and investigated.

The geometry of finite dimensional real Banach (or normed) spaces,  
also called Minkowski geometry, is studied in the basic references \cite{BMS}, \cite{JMR}, \cite{MS}, \cite{MSW}, \cite{MW} and \cite{T}, and recently it has strong relations to fields like optimization, discrete and computational geometry, 
convexity, convex and functional analysis, approximation theory and so on.
Also from the viewpoint of differential geometry it is natural to  
develop geometric concepts for norms. Regarding curvature types of planar curves this was, as already mentioned, done in \cite{BMS}, and for surfaces in  
higher dimensions we refer to \cite{BMT} and \cite{JMR}. Deleting the symmetry axiom  
(i.e., going from norms to gauges), one can find almost no analogous  
references. We could locate the single paper \cite{G} which contains  
various results on curves and hypersurfaces derived explicitly for gauge planes and  
spaces in the spirit of differential geometry, and in \cite{JMR} certain  
multifocal hypersurfaces (e.g., polyellipsoids) are studied for norms and gauges.

A function $F: {\mathbb R}^n\rightarrow {\mathbb R}$ on the $n$-dimensional linear space ${\mathbb R}^n$ is called a \emph{convex distance function}, or a \emph{gauge}, if it satisfies the following conditions:

(i) $F(x)\geq 0$ for $x \in {\mathbb R}^n$, and $F(x) = 0 \Leftrightarrow x = 0$,

(ii) $F(\lambda x) = \lambda F(x)$ for $x \in {\mathbb R}^n, \lambda > 0$,

(iii) $F(x+y) \leq F(x)+F(y)$ for $x, y \in {\mathbb R}^n$. \\
Having this notion, the space $({\mathbb R}^n, F)$ is called an $n$-dimensional \emph{generalized Minkowski space} or \emph{gauge space} (cf. \cite{JMR}), hence being a direct extension of the concept of 
normed or Minkowski space (or linear Finsler space). We expect that the study of generalized Minkowski spaces will give new insights important also for Finsler geometry (cf. \cite{BCS} and \cite{CS}). Clearly, for $n = 2$ 
\emph{generalized Minkowski planes} or \emph{gauge planes} are obtained.

Our paper is organized as follows. In Section 2 we give some basic  
notions for and facts on generalized Minkowski planes. In particular,  
the notions of associated gauge and Birkhoff orthogonality are  
introduced. In Section 3 we directly generalize the main results from \cite{BMS}, namely by deriving the four announced curvature types of curves in gauge planes. Also relations between these four curvature  
types are presented there. Applying then our concepts, we prove in Section 4  
the bi-directional relation between evolutes and involutes of given  
curves in gauge planes. All these results differ from those published  
on curvatures, evolutes, and involutes in \cite{G}, see the following Remark.\\

\textbf{Remark.} In \cite{G} Guggenheimer discussed two curvature types of curves in generalized Minkowski planes which correspond to the Minkowski and the normal curvature in our setting. But different to his approach, our formulation is based on the notion of associated gauge, and explicit computational formulas for curvatures are obtained. Also, Guggenheimer defined a type of evolutes which is different from ours. For his setting, only a one-sided relation between evolutes and involutes can be shown.

\section{Basic facts}

Let $({\mathbb R^2}, F)$ be a generalized Minkowski plane whose \emph{unit disk} $B$ and whose \emph{unit circle} $S$ are defined by
\[B = \{x \in {\mathbb R^2}; F(x) \leq 1\}, \ \ \ \ S = \{x \in {\mathbb R^2}; F(x) = 1\}\,; \]
here $B$ is a compact, convex set having the origin $0$ as interior point and $S$ as its boundary. This is equivalent to the property that the considered plane is equipped with 
a convex distance function $F$ as defined in the introduction.

Let $[\cdot, \cdot ]: {\mathbb R^2}\times {\mathbb R^2} \rightarrow {\mathbb R}$ be a map given by
\[[x, y] = \left| \begin{array} {cc} x_1 & y_1 \\ x_2 & y_2 \end{array} \right|, \]
where $x = (x_1, x_2)$, $y = (y_1, y_2)$ and $|\ \ |$ denotes the usual determinant. Based on this, let
\[F_a(x) = \mbox{sup}\{[y, x]; F(y) = 1\}, \]
which is also a gauge. We call $F_a$ the \emph{gauge associated to} $F$. Then $[y, x] \leq F(y)F_a(x)$. Considering the orientation, we say that $x$ is \emph{Birkhoff orthogonal} to $y$ (denoted by $x \dashv_{B} y$) if $[x, y] = F(x)F_a(y)$ (cf. \cite{MS} for the subcase of normed planes).

We set
\[S_a = \{x \in {\mathbb R^2}; F_a(x) = 1\} \]
and let
\[F_{a,a}(x) = (F_a)_{a}(x) = \mbox{sup}\{[y, x] ; F_a(y) = 1\}. \]

\vspace{2mm}

{\bf Example 1.} We consider a Randers norm on ${\mathbb R}^2$ given by
\[F(x) = F(x_1, x_2) = \sqrt{x_1^2+x_2^2}+bx_1, \]
where $|b| < 1$. The equation $F(y) = 1$ is equivalent to
\[\{(1-b^2)y_1+b\}^2+(1-b^2)y_2^2 = 1. \]
We can compute that
\[F_a(x) = \frac{1}{1-b^2}\sqrt{(1-b^2)x_1^2+x_2^2}-\frac{b}{1-b^2}x_2, \]
which is also a Randers norm. The equation $F_a(y) = 1$ is equivalent to
\[y_1^2+(y_2-b)^2 = 1, \]
and we can show that
\[F_{a,a}(x) = \sqrt{x_1^2+x_2^2}-bx_1 = F(-x). \]

\vspace{2mm}

{\bf Example 2.} Let $D$ be a triangle in ${\mathbb R^2}$ containing the origin as interior point and having the vertices $(p_1, p_2)$, $(q_1, q_2)$ and $(r_1, r_2)$. Then we can interpret $D$ as unit disk of a gauge $F$. We can see that the unit disk $B_a$ for $F_a$ is represented by
\[p_1 x_2-p_2 x_1 \leq 1,\ \ \ \ q_1 x_2-q_2 x_1 \leq 1,\ \ \ \ r_1 x_2-r_2 x_1 \leq 1, \]
with vertices
\[\frac{1}{p_1 q_2-p_2 q_1}(q_1-p_1, q_2-p_2), \ \ \ \ \frac{1}{q_1 r_2-q_2 r_1}(r_1-q_1, r_2-q_2), \]
\[\frac{1}{r_1 p_2-r_2 p_1}(p_1-r_1, p_2-r_2). \]
Similarly computing, we find that the unit disk $B_{a,a}$ for $F_{a,a}$ is the triangle with vertices $-(p_1, p_2)$, $-(q_1, q_2)$ and $-(r_1, r_2)$. This implies that $F_{a,a}(x) = F(-x)$. \\

In general we get the following

\begin{prop} Let $F$ be a gauge on ${\mathbb R^2}$. Then we have:

(i) $F_{a,a}(x) = F(-x)$ for any $x \in {\mathbb R^2}$.

(ii) $x \dashv_{B} y$ is equivalent to the fact that $y$ is Birkhoff orthogonal to $-x$ with respect to $F_a$ (denoted by $y \dashv_{B}^{a}(-x)$).
\end{prop}

Proof. (i) Since
\[[y, x] = [-x, y] \leq F(-x)F_{a}(y), \]
we have $F_{a,a}(x) \leq F(-x)$. Conversely, for any $x$ we can choose $\tilde{y} \neq 0$ so that $[-x, \tilde{y}] = F(-x)F_{a}(\tilde{y})$. Then
\[F_{a,a}(x) \geq \left[\frac{\tilde{y}}{F_{a}(\tilde{y})}, x \right] = \frac{1}{F_{a}(\tilde{y})}[-x, \tilde{y}] = F(-x). \]
Thus we have $F_{a,a}(x) = F(-x)$.

(ii) By (i), the Birkhoff orthogonality condition $[x, y] = F(x)F_{a}(y)$ is rewritten as $[y, -x] = F_{a}(y)F_{a,a}(-x)$. So $x \dashv_{B} y$ is equivalent to the fact that $y$ is Birkhoff orthogonal to $-x$ with respect to $F_a$.
\hfill $\Box$ \\

\textbf{Remark.} As in Example 2, for a convex polygon containing the origin as interior point we can show that $F_{a,a}(x) = F(-x)$. Thus, suitable approximation of general convex compact domains containing the origin as interior point by convex polygons yields that we have $F_{a,a}(x) = F(-x)$ also for this general setting. This gives another proof of Prop. 2.1(i). \\

\textbf{Remark.} When $F$ is symmetric (that is, $F(x) = F(-x)$), then we have $F_{a,a}(x) = F(x)$ (cf. \cite{MS}).

\section{Curvature types}

Let $({\mathbb R^2}, F)$ be a generalized Minkowski plane. In the following, we assume that its unit circle $S$ is smooth and strictly convex. Let $\gamma(s)$ be an oriented smooth curve in $({\mathbb R^2}, F)$ with arc-length parameter $s$ such that $F(\gamma'(s)) = 1$.

Let $\varphi(t)$ be a counter-clockwise parametrization of $S$ such that $[\varphi(t), \varphi'(t)] > 0$. Set
\[u(t) = \int_{t_0}^t [\varphi(\tau), \varphi'(\tau)]d\tau, \]
and let $\varphi(u)$ be a parametrization of $S$ by $u$. Then we have
\[\varphi(u) \dashv_{B} \varphi'(u), \ \ \ \ F_a(\varphi'(u)) = 1. \]
We can write $\gamma'(s) = \varphi(u(s))$, and we define the \emph{Minkowski curvature} by $k_m(s) := u'(s)$. We set $n_{\gamma}(s) = \varphi'(u(s)) \in S_a$, which we call the \emph{right normal vector field}. Then we have
\[\gamma'(s) \dashv_{B} n_{\gamma}(s), \ \ \ \ [\gamma'(s), n_{\gamma}(s)] = 1, \ \ \ \ \gamma''(s) = k_m(s)n_{\gamma}(s). \]

Let $\psi(v)$ be a counter-clockwise parametrization of the unit circle $S_a$ for $F_a$ such that $F_{a,a}(\psi'(v)) = F(-\psi'(v)) = 1$. Then we can write $n_{\gamma}(s) = \psi(v(s))$. We define the \emph{normal curvature} by $k_n(s) := v'(s)$.

\begin{lem}
$\gamma'(s) = -\psi'(v(s))$.
\end{lem}

Proof. By the condition, we have $\psi(v)\dashv_{B}^{a} \psi'(v)$. By (ii) of Prop. 2.1, this is equivalent to $-\psi'(v) \dashv_{B} \psi(v)$. Thus we have
\[\gamma'(s) \dashv_{B} n_{\gamma}(s),\ \ \ \ -\psi'(v(s))\dashv_{B} n_{\gamma}(s). \]
Since $F(\gamma'(s)) = F(-\psi'(v(s))) = 1$ and $S$ is strictly convex, we get $\gamma'(s) = -\psi'(v(s))$. \hfill $\Box$ \\

By Lemma 3.1, we have
\[n_{\gamma}'(s) = -k_n(s) \gamma'(s). \]

\vspace{3mm}

Let $\varphi(t)$ be a counter-clockwise parametrization of $S$ by the arc-length parameter $t$ such that $F(\varphi'(t)) = 1$. Then we can write $\gamma'(s) = \varphi'(t(s))$, and we define the \emph{circular curvature} by $k_c(s) := t'(s)$. We call $\varphi(t(s)) \in S$ the \emph{left normal vector field}, which satisfies $[\gamma'(s), \varphi(t(s))] < 0$. Now we define the \emph{evolute} $E$ of $\gamma$ by
\[E(s) = \gamma(s)-\frac{1}{k_c(s)}\varphi(t(s)). \]
Then
\[E'(s) = -\left(\frac{1}{k_c(s)}\right)' \varphi(t(s)). \]
So $E$ is the \emph{envelop} of left normal lines.

For the counter-clockwise parametrization $\varphi(t)$ of $S$ by the arc-length parameter $t$ with $F(\varphi'(t)) = 1$, we can also write $\gamma'(s) = \varphi(\hat{t}(s))$, and we define the \emph{arc-length curvature} by $k_{l}(s) := \hat{t}'(s)$. Thus we can define the \emph{involute} $I$ of $\gamma$ by
\[I(s) = \gamma(s)+(c-s)\gamma'(s) = \gamma(s)+(c-s)\varphi(\hat{t}(s)). \]

Next we discuss relations between the four obtained types of curvatures.

\begin{prop}
$k_{l}(s) = F(n_{\gamma}(s)) k_m(s)$. 
\end{prop}

Proof. Since 
\[\gamma''(s) = k_m(s)n_{\gamma}(s) = k_{l}(s) \varphi'(\hat{t}(s)), \]
\[\gamma'(s) \dashv_{B} n_{\gamma}(s),\ \ \ \ \gamma'(s) = \varphi(\hat{t}(s)) \dashv_{B} \varphi'(\hat{t}(s)), \]
\[n_{\gamma}(s) \in S_a, \ \ \ \ \varphi'(\hat{t}(s)) \in S, \]
we have 
\[\varphi'(\hat{t}(s)) = \frac{n_{\gamma}(s)}{F(n_{\gamma}(s))}, \ \ \ \ k_{l}(s) = F(n_{\gamma}(s)) k_m(s). \qquad \qquad \Box \]

\begin{prop}
The circular curvature $k_c$ is the normal curvature $k_{n}^{a}$ with respect to the associated gauge $F_a$.
\end{prop}

Proof. Let $s_a$ be the arc length of $\gamma$ with respect to $F_a$. Let $n_{\gamma}^a$ be the right normal vector field to $\gamma$ with respect to $F_a$. Since $n_{\gamma}^a \in S_{a,a} = (S_a)_a$, we have $F_{a,a}(n_{\gamma}^a) = F(-n_{\gamma}^a) = 1$. Then
\[\frac{d\gamma}{ds_a}(s_a) \dashv_B^a n_{\gamma}^a (s_a), \]
and $\gamma'(s) \dashv_B^a n_{\gamma}^a (s_a(s))$. By (ii) of Prop. 2.1 and $\gamma'(s) = \varphi'(t(s))$, we have
\[-n_{\gamma}^a (s_a(s)) \dashv_B \varphi'(t(s)). \]
On the other hand,
\[\varphi(t(s)) \dashv_B \varphi'(t(s)). \]
Since $-n_{\gamma}^a (s_a(s)) \in S$, $\varphi(t(s)) \in S$ and $S$ is strictly convex, we find that
\[-n_{\gamma}^a (s_a(s)) = \varphi(t(s)). \]
Differentiation with respect to $s$ yields
\[k_n^a(s_a(s)) \frac{d\gamma}{ds_a}(s_a(s)) \cdot \frac{ds_a}{ds} = \varphi'(t(s)) t'(s) \]
and
\[k_n^a(s_a(s)) \gamma'(s) = k_c(s) \gamma'(s). \]
Thus we get $k_c(s) = k_n^a(s_a(s))$. \hfill $\Box$ \\

We assume that $F$ is smooth on ${\mathbb R}^2 \setminus \{0\}$ and give formulas for curvatures. Let $(x_1, x_2)$ be the standard coordinates on $({\mathbb R}^2, F)$. From the positive homogeneity of $F$, we have 

(a) \ \ $\displaystyle{ \sum_{i=1}^{2}t_i F_{x_i}(t_1, t_2) = F(t_1, t_2) }$ \\
and

(b) \ \ $\displaystyle{ \sum_{i=1}^{2} t_i F_{x_i x_j}(t_1, t_2) = 0 }$ \\
for $(t_1, t_2) \in {\mathbb R}^2 \setminus \{0\}$ (cf. Chap. 1 of \cite{BCS}). 

Let $\gamma(\tau) = (\gamma_1(\tau), \gamma_2(\tau))$ be an oriented smooth curve in $({\mathbb R}^2, F)$ with arbitrary parameter $\tau$. Let $s$ be the arc-length parameter. Then 
\[\frac{d\gamma}{ds} = \frac{1}{F(\gamma_1', \gamma_2')} (\gamma_1', \gamma_2'). \]
We note that the vector 
\[(-F_{x_2}(\gamma_1', \gamma_2'), F_{x_1}(\gamma_1', \gamma_2')) \]
is parallel to the tangential line of $S$ at $d\gamma/ds$. By the definition of the associated gauge $F_a$, we have 
\[F_{a}(-F_{x_2}(\gamma_1', \gamma_2'), F_{x_1}(\gamma_1', \gamma_2')) = 1, \]
where we use (a). So the right normal vector field $n_{\gamma}$ is given by 
\[n_{\gamma} = (-F_{x_2}(\gamma_1', \gamma_2'), F_{x_1}(\gamma_1', \gamma_2')). \] 

Using (a), we can compute 
\[\frac{d^2\gamma}{ds^2} = \frac{1}{(F(\gamma_1', \gamma_2'))^3} (\gamma_{1}'' \{F(\gamma_1', \gamma_2')-\gamma_1' F_{x_1}(\gamma_1', \gamma_2')\}-\gamma_1' \gamma_2'' F_{x_2}(\gamma_1', \gamma_2'),  \]
\[ \hspace{2.5cm} \gamma_{2}'' \{F(\gamma_1', \gamma_2')-\gamma_2' F_{x_2}(\gamma_1', \gamma_2')\}-\gamma_1'' \gamma_2' F_{x_1}(\gamma_1', \gamma_2') ) \]
\[= \frac{\gamma_{1}' \gamma_{2}''-\gamma_{1}'' \gamma_{2}'}{(F(\gamma_1', \gamma_2'))^3} (-F_{x_2}(\gamma_1', \gamma_2'), F_{x_1}(\gamma_1', \gamma_2')) = \frac{\gamma_{1}' \gamma_{2}''-\gamma_{1}'' \gamma_{2}'}{(F(\gamma_1', \gamma_2'))^3}\cdot n_{\gamma}. \]
Thus we get

\begin{thm}
The Minkowski curvature is computed by 
\[k_m = \frac{\gamma_{1}' \gamma_{2}''-\gamma_{1}'' \gamma_{2}'}{(F(\gamma_1', \gamma_2'))^3}. \]
\end{thm}

\vspace{3mm}

And by Proposition 3.1 and Theorem 3.1, we obtain

\begin{cor}
The arc-length curvature is given by 
\[k_l = \frac{\gamma_{1}' \gamma_{2}''-\gamma_{1}'' \gamma_{2}'}{(F(\gamma_1', \gamma_2'))^3} \cdot F(-F_{x_2}(\gamma_1', \gamma_2'), F_{x_1}(\gamma_1', \gamma_2')). \]
\end{cor}

\vspace{3mm}

Next we have 
\[\frac{dn_{\gamma}}{ds} = \frac{1}{F(\gamma_1', \gamma_2')} (-(F_{x_2}(\gamma_1', \gamma_2'))', (F_{x_1}(\gamma_1', \gamma_2'))' ). \]
Differentiating the relation 
\[\gamma_1' F_{x_1}(\gamma_1', \gamma_2')+\gamma_2' F_{x_2}(\gamma_1', \gamma_2') = F(\gamma_1', \gamma_2') \]
from (a), we have 
\[\gamma_1' (F_{x_1}(\gamma_1', \gamma_2'))'+\gamma_2' (F_{x_2}(\gamma_1', \gamma_2'))' = 0. \]
Therefore $\gamma_{1}' \neq 0$ yields
\[\frac{dn_{\gamma}}{ds} = -\frac{(F_{x_2}(\gamma_1', \gamma_2'))'} {\gamma_1' F(\gamma_1', \gamma_2')} (\gamma_1', \gamma_2') = -\frac{(F_{x_2}(\gamma_1', \gamma_2'))'} {\gamma_1'} \cdot \frac{d\gamma}{ds}. \]
Similarly, when $\gamma_{2}' \neq 0$, then we obtain
\[\frac{dn_{\gamma}}{ds} = \frac{(F_{x_1}(\gamma_1', \gamma_2'))'}{\gamma_2'} \cdot  \frac{d\gamma}{ds}. \]
Thus we have the following theorem.

\begin{thm}
(i) If $\gamma_{1}' \neq 0$, then we get
\[k_n = \frac{(F_{x_2}(\gamma_1', \gamma_2'))'}{\gamma_1'}. \]
(ii) If $\gamma_{2}' \neq 0$, then we have
\[k_n = -\frac{(F_{x_1}(\gamma_1', \gamma_2'))'}{\gamma_2'}. \]
\end{thm}

\vspace{3mm}

\textbf{Remark.} The above computation on $k_n$ is related to (b). \\

By Proposition 3.2 and Theorem 3.2, the following corollary is obtained.

\begin{cor}
(i) If $\gamma_{1}' \neq 0$, then we obtain
\[k_c = \frac{((F_a)_{x_2}(\gamma_1', \gamma_2'))'}{\gamma_1'}. \]
(ii) If $\gamma_{2}' \neq 0$, then we see that
\[k_c = -\frac{((F_a)_{x_1}(\gamma_1', \gamma_2'))'}{\gamma_2'}. \]
\end{cor}

\vspace{3mm}

\textbf{Remark.} In general, the formula for $k_c$ is not explicit because $F_a$ is used. \\

\textbf{Remark.} (i) In the Euclidean case, those formulas reduce to the classical one. 

(ii) The analogous discussion in normed planes can be seen in Section 5 of \cite{BMS}, where radial coordinates are used.

\section{Evolutes and involutes}

Let $\gamma(s)$ be an oriented smooth curve in $({\mathbb R^2}, F)$ with arc-length parameter $s$ so that $F(\gamma'(s)) = 1$, as in Section 3. For an oriented curve $\alpha(t)$, we define the \emph{reverse curve} $\alpha^{-}$ by $\alpha^{-}(t) = \alpha(-t)$.

\begin{thm}
(i) If $(c-s)k_{l}(s) > 0$, then the evolute of the involute $I$ of $\gamma$ is $\gamma$.

(ii) If $(c-s)k_{l}(s) < 0$, then the reverse evolute of the reverse involute $I^{-}$ of $\gamma$ is $\gamma$.
\end{thm}

Proof. (i) First we have
\[I'(s) = (c-s)k_{l}(s)\varphi'(\hat{t}(s)). \]
Let $s^{\ast}$ denote the arc length of $I$. Since $(c-s)k_{l}(s) > 0$, we have
\[\frac{ds^{\ast}}{ds} = F(I'(s)) = (c-s)k_{l}(s),\ \ \ \ \frac{dI}{ds^{\ast}} = \varphi'(\hat{t}(s)). \]
The circular curvature $k^{\ast}_c$ of $I$ is given by
\[k^{\ast}_{c} = \frac{d\hat{t}}{ds^{\ast}} = \frac{d\hat{t}}{ds}\frac{ds}{ds^{\ast}} = \frac{1}{c-s}. \]
So the evolute $E_I$ of $I$ satisfies
\[E_I(s) = I(s)-\frac{1}{k^{\ast}_{c}(s)}\varphi(\hat{t}(s)) =\gamma(s). \]

(ii) Since $(c-s)k_{l}(s) < 0$ and letting $s = -\sigma$, we have $(c+\sigma)k_{l}(-\sigma) < 0$. The reverse involute $I^{-}$ is given by
\[I^{-}(\sigma) = I(-\sigma) = \gamma(-\sigma)+(c+\sigma)\varphi(\hat{t}(-\sigma)). \]
Then
\[(I^{-})'(\sigma) = -(c+\sigma)k_{l}(-\sigma)\varphi'(\hat{t}(-\sigma)). \]
Let $s^{\ast}$ denote the arc length of $I^{-}$. Then
\[\frac{ds^{\ast}}{d\sigma} = F((I^{-})'(\sigma)) = -(c+\sigma)k_{l}(-\sigma),\ \ \ \ \frac{dI^{-}}{ds^{\ast}} = \varphi'(\hat{t}(-\sigma)). \]
The circular curvature $k^{\ast}_c$ of $I^{-}$ is given by
\[k^{\ast}_{c} = \frac{d}{ds^{\ast}}(\hat{t}(-\sigma)) = \frac{d}{d\sigma}(\hat{t}(-\sigma)) \frac{d\sigma}{ds^{\ast}} = \frac{1}{c+\sigma}. \]
So the evolute $E_{(I^{-})}$ of $I^{-}$ satisfies
\[E_{(I^{-})}(\sigma) = I^{-}(\sigma)-\frac{1}{k^{\ast}_{c}(\sigma)}\varphi(\hat{t}(-\sigma)) = \gamma(-\sigma). \]
Now, by $\sigma = -s$ (that is, reversing the orientation), we get
\[E_{(I^{-})}(-s) = \gamma(s). \qquad \qquad \qquad \qquad \qquad \qquad \Box
\]

\begin{thm}
(i) If $k_c'(s) > 0$, then an involute of the evolute $E$ of $\gamma$ is $\gamma$.

(ii) If $k_c'(s) < 0$, then a reverse involute of the reverse evolute $E^{-}$ of $\gamma$ is $\gamma$.
\end{thm}

Proof.
(i) We have
\[E'(s) = -\left(\frac{1}{k_c(s)}\right)' \varphi(t(s)) =  \frac{k_{c}'(s)}{(k_c(s))^2} \varphi(t(s)). \]
Let $s^{\ast}$ denote the arc length of $E$. Since $k_{c}'(s) > 0$, we have
\[\frac{ds^{\ast}}{ds} = F(E'(s)) = -\left(\frac{1}{k_c(s)}\right)' = \frac{k_{c}'(s)}{(k_c(s))^2}, \]
\[\frac{dE}{ds^{\ast}} = \varphi(t(s)), \]
and
\[s^{\ast} = -\int\left(\frac{1}{k_c(s)}\right)' ds = -\frac{1}{k_c(s)}+c. \]
So an involute $I_E$ of $E$ satisfies
\[I_E(s) = E(s)+(c-s^{\ast})\frac{dE}{ds^{\ast}} = \gamma(s). \]

(ii) Since $k_c'(s) < 0$ and letting $s = -\sigma$, we have $k_c'(-\sigma) < 0$. The reverse evolute $E^{-}$ is given by
\[E^{-}(\sigma) = E(-\sigma) = \gamma(-\sigma)-\frac{1}{k_c(-\sigma)}\varphi(t(-\sigma)). \]
Then
\[(E^{-})'(\sigma) = -\left(\frac{1}{k_c(-\sigma)}\right)' \varphi(t(-\sigma)) =  -\frac{k_{c}'(-\sigma)}{(k_c(-\sigma))^2} \varphi(t(-\sigma)). \]
Let $s^{\ast}$ denote the arc length of $E^{-}$. Then
\[\frac{ds^{\ast}}{d\sigma} = F((E^{-})'(\sigma)) = -\left(\frac{1}{k_c(-\sigma)}\right)' = -\frac{k_{c}'(-\sigma)}{(k_c(-\sigma))^2}, \]
\[\frac{dE^{-}}{ds^{\ast}} = \varphi(t(-\sigma)), \]
and
\[s^{\ast} = -\int\left(\frac{1}{k_c(-\sigma)}\right)' d\sigma = -\frac{1}{k_c(-\sigma)}+c. \]
So an involute $I_{(E^{-})}$ of $E^{-}$ satisfies
\[I_{(E^{-})}(\sigma) = E^{-}(\sigma)+(c-s^{\ast})\frac{dE^{-}}{ds^{\ast}} = \gamma(-\sigma). \]
Reversing the orientation by $\sigma = -s$, we get
\[I_{(E^{-})}(-s) = \gamma(s). \qquad \qquad \qquad \qquad \qquad \qquad \Box
\]


{\small
\begin{thebibliography}{99}
\bibitem{BMS} V. Balestro, H. Martini and E. Shonoda, Concepts of curvatures in normed planes, \emph{Expo. Math.}, to appear, arXiv:1702.01449.
\bibitem{BMT} V. Balestro, H. Martini and R. Teixeira, Differential geometry of immersed surfaces in three-dimensional normed spaces, arXiv:1707.04226.
\bibitem{BCS} D. Bao, S. S. Chern and Z. Shen, \emph{An Introduction to Riemann-Finsler Geometry}, Springer, 2000.
\bibitem{CS} S. S. Chern and Z. Shen, \emph{Riemann-Finsler Geometry}, World Sci., 2005.
\bibitem{G} H. Guggenheimer, Pseudo-Minkowski differential geometry, \emph{Ann. Mat. Pure Appl.} \textbf{70} (1965), 305-370.
\bibitem{JMR} T. Jahn, H. Martini and C. Richter, Bi- and multifocal curves and surfaces for gauges, \emph{J. Convex Anal.} \textbf{23} (2016), 733-774.
\bibitem{MS} H. Martini and K. J. Swanepoel, Antinorms and Radon curves, \emph{Aequationes Math.} \textbf{72} (2006), 110-138.
\bibitem{MSW} H. Martini, K. J. Swanepoel and G. Weiss, The geometry of 
Minkowski spaces - a survey. I. \emph{Expo. Math.} \textbf{19} (2001), 
no. 2, 97-142.
\bibitem{MW} H. Martini and S. Wu, Classical curve theory in normed planes, \emph{Comput. Aided Geom. Design} \textbf{31} (2014), 373-397.
\bibitem{T} A. C. Thompson, \emph{Minkowski Geometry}, Cambridge University  
Press, Cambridge, 1996.
\end{thebibliography}
}
\end{document}